\documentclass[11pt,a4]{article}
\usepackage{threeparttable}
\usepackage{amsmath,amssymb,amsthm}
\usepackage{color}
\usepackage{collcell}
\usepackage[english]{babel}
\usepackage{graphicx}
\usepackage{natbib}
\usepackage{subfig}
\usepackage{threeparttable}
\usepackage{multirow}

\usepackage{booktabs}

\newtheorem{definition}{Definition}%

\newcommand{\de}{\mathrm d}

\newcommand{\DPP}{\textsc{dpp}}
\newcommand{\stDPP}{\textsc{stdpp}}

\newcommand{\SOIRS}{\textsc{soirs}}
\usepackage{mathtools}

\usepackage{enumerate}

\usepackage[utf8]{inputenc}
\usepackage{fourier} 
\usepackage{array}
\usepackage{makecell}

\usepackage[utf8]{inputenc}
\usepackage{rotating}
\usepackage{threeparttable}
\usepackage{amsmath,amssymb,amsthm}
\usepackage{threeparttable}
\usepackage{amsmath,amssymb,amsthm}
\usepackage{color}
\usepackage{collcell}
\usepackage[english]{babel}
\usepackage{graphicx}
\usepackage{natbib}
\usepackage{subfig}
\usepackage{threeparttable}
\usepackage{multirow}
\usepackage{booktabs}

\usepackage{setspace}

\usepackage{geometry}
\newtheorem{def.}{Definition}[section]
\newtheorem{examp.}{Example}[section]

\RequirePackage[colorlinks,citecolor=blue,urlcolor=blue]{hyperref} 

\begin{document}

\title{
 Spatio-temporal determinantal point processes 
}
\date{}

\author{}
\maketitle
  
\begin{center}
  {\small\bf Nafiseh Vafaei}\\
 \thanks{Department of Computer and Statistics Sciences, Faculty of Sciences, University of Mohaghegh Ardabili, Ardabil, Iran\\E-mail: N.Vafaei@uma.ac.ir}\\
 
  {\small\bf Mohammad Ghorbani}\footnote{Corresponding Author}\\\thanks{Department of Engineering Sciences and Mathematics, Lule\r{a} University of Technology, Sweden\\E-mail: 
    mohammad.ghorbani@ltu.se} \\

   {\small\bf Masoud Ganji}\\
 \thanks{Department of Computer and Statistics Sciences, Faculty of Sciences, University of Mohaghegh Ardabili, Ardabil, Iran
 \\E-mail: mganji@uma.ac.ir}\\ 

  {\small\bf Mari Myllym\"{a}ki}\\
 \thanks{Natural Resources Institute Finland (Luke),  Helsinki, Finland.\\E-mail: mari.myllymaki@luke.fi}\\ 

\end{center}
\begin{abstract}
Determinantal point processes are models for regular spatial point patterns, with appealing probabilistic properties. 
We present their spatio-temporal counterparts and give examples of these models, based on spatio-temporal covariance functions which are separable and non-separable in space and time.
\end{abstract}



\section{Introduction}\label{sec:ASS}
Spatio-temporal point processes are random countable subsets $X$ of $\mathbb R^2\times\mathbb R^+$, where a point $(u,t)\in X$ corresponds to an event $u\in\mathbb R^2$ occurring at time $t\in\mathbb R^+$. We assume that the points do not overlap, i.e.\ $(u_1, t_1)\ne (u_2, t_2)$. 
Examples of such events are the occurrence of epidemic diseases (such as corona or flu), sightings or births of a species, the occurrence of fires, earthquakes, tsunamis, or volcanic eruptions.
We are interested here in spatio-temporal regular point processes, where neighbouring points in the process tend to repel each other.

In the context of spatial point processes, Gibbs point processes including Markov point processes and pairwise interaction point process models 
are generally used to model repulsiveness. 
Another class of regular spatial point process models are determinantal point processes (\DPP{}s), which have their origin in quantum physics. They were first identified as a class by \cite{Macchi:75}, who called them {\em fermion processes} because they reflect the distributions of fermion systems in thermal equilibrium that exhibit repulsive behaviour. 
\DPP{}s have been extensively studied in probability theory 
and have found applications in random matrix theory, 
quantum physics, 
wireless network modelling, 
Monte Carlo integration, 
and machine learning \citep[see e.g.][]{Lavancier:etal:15}. 
Recently, \cite{Lavancier:etal:15} studied the statistical properties of \DPP{}s.

Some regular spatial point process models have already their spatio-temporal counterparts, but 
likelihood-based inference or simulation of these models is usually complicated and time-consuming.
To circumvent these challenges, our objective here is to introduce the spatio-temporal determinantal point processes (\stDPP{}s) and study their properties, which is an open problem according to \citet[ pages 875-876]{Lavancier:etal:15}.
We present the basic properties of these processes and give examples of them based on separable and non-separable spatio-temporal covariance functions. We derive the key summary characteristics for these examples that can be used, for example, for model fitting and evaluation. 
\section{Basic concepts and statistical properties}\label{sec:stDPP}
Assume that $X$ is a spatio-temporal point process with $n$th-order product density $\rho^{(n)}$,   $n\ge 1$ which describes the frequency of possible configurations of $n$ points. Suppose that $B_1,\ldots,B_n$ are pairwise disjoint cylindrical regions having infinitesimal volumes $dV_1, \ldots,dV_n$ and containing the points $((u_1,t_1), \dots, (u_n, t_n))$, respectively.
Then, $\rho^{(n)}\big((u_1,t_1),\ldots,(u_n,t_n)\big)dV_1,\ldots dV_n$ is the probability that $X$ has a point in each of  $B_1,\ldots,B_n$. 

\begin{definition}\label{def1}
Let $C$ be a kernel function from $(\mathbb{R}^2\times\mathbb{R}^+)\times(\mathbb{R}^2\times\mathbb{R}^+)$ to  $\mathbb{R}$.
We say that a spatio-temporal point process $X$ is a determinantal point process with kernel $C$ and write $X\sim \stDPP{}(C)$, if its $n$th-order product density function is given by
\begin{align*}
\rho^{(n)}((u_1, t_1),\ldots, (u_n,t_n)) =
\det\{C((u_i,t_i), (u_j, t_j))\}_{1\leq i,j \leq n},
\end{align*}
for $(u_i,t_i)\in (\mathbb{R}^{2}\times\mathbb{R}^+)$  and $n=1,2,\ldots$, where $\{C((u_i,t_i),(u_j,t_j)\}_{1\leq i,j\leq n}$ is the $n\times n$ matrix with $C((u_i,t_i), (u_j, t_j))$  as its $(i,j)$-th entry and $\det \{A\}$ is the determinant of the   matrix $A$. 
\end{definition}
The point process $X$ is well-defined if, for each $n$, $\rho^{(n)}((u_1,t_1),\cdots,(u_n,t_n))\geq 0$ for all $\{(u_i,t_i)\}_{i=1}^n$. This implies that, in Definition~\ref{def1}, 
$\mathbf{C}= \{C((u_i,t_i),(u_j,t_j)\}_{1\leq i,j\leq n}$ should be a non-negative definite matrix.
Then all eigenvalues of the matrix $\mathbf{C}$ are non-negative and thus its determinant  is also non-negative (this follows from the fact that the determinant of a matrix is equal to the product of its eigenvalues) \citep{Rao:Rao:98}. Therefore, covariance functions are possible choices for the kernel $C$.
Denoting the eigenvalues of $\mathbf{C}$ by $\lambda_l$ $(l=1, 2, \ldots)$, another condition for the existence of \stDPP{}($C$) is that $\lambda_l\le 1$ $(l=1, 2, \ldots)$. This is a straightforward extension from the spatial setting, see details in \cite{Lavancier:etal:15} and the references therein. 
The Poisson process with intensity $\rho(u,s)$ results as a special case of \stDPP{}($C$) with setting
$C((u,s),(u,s)) = \rho(u,s)$ and $C((u,s),(v,t)) = 0 \,\, \mbox{if}\,\, u\ne v\, \mbox{or}\, s\ne t$.

By Definition~\ref{def1}, the moment characteristics of arbitrary order $n$ $(n\ge 1)$ can be easily attained for \stDPP{}s, e.g. the intensity function is given by $\rho(u,t)=C((u,t), (u,t))$ 
and the second-order product density is given by
\begin{align*}
   \rho^{(2)}((u_1, t_1), (u_2,t_2)) &=C\left((u_1, t_1), (u_1,t_1)\right)C((u_2, t_2), (u_2,t_2))\cr &-C((u_1, t_1), (u_2,t_2))C((u_1, t_1), (u_2,t_2)). 
\end{align*}
In what follows, we assume that the kernel function is of the form
\begin{align*} 
C((u_1,t_1),(u_2,t_2))=\sqrt{\rho(u_1,t_1)} R\Big((u_1,t_1) /\alpha_s,(u_2,t_2)/\alpha_t\Big)\sqrt{\rho(u_2,t_2)},
\end{align*}
where $\rho$ plays the role of the intensity function of the process, 
$\alpha_s>0$ and $\alpha_t >0$ are the spatial and  temporal correlation parameters respectively,
and $R(\cdot)$ is the correlation function correspondent to $C$. Then,  
the (inhomogeneous space-time) pair correlation function  is   given by 
\begin{align}\label{eq:gfun}
g((u_1,t_1), (u_2,t_2))&=\frac{ \rho^{(2)}((u_1, t_1), (u_2,t_2))}{\rho(u_1, t_1)\rho(u_2, t_2)}\cr
&=1- \Big \lvert R\Big((u_1,t_1) /\alpha_s,(u_2,t_2)/\alpha_t\Big)\Big \rvert^2\le 1. 
\end{align}
Since $g\le 1$, the points of $X$ repel each other, which is a characteristic of regular point patterns.

In general, a \stDPP{} is called second-order intensity reweighted stationary (\SOIRS{})  if its space-time pair correlation function is a function of the spatial difference $u=u_2-u_1$ and the temporal difference $t=t_2-t_1$ \citep[see, e.g.][and the references therein]{ghorbani:13}. Thus, a \stDPP{} with kernel function $C$ is \SOIRS{} if the correlation function $R$ is a function of $u$ and $r$ only, i.e.\ if $R$ takes the form $R((u_1,t_1), (u_2,t_2))=R_0(u, t)$ for some function $R_0$. 
For a \SOIRS{} \stDPP{}(C), $R_0(0,0)=1$, and hence $\rho(u,t)=C(0,0)$. 
If in addition, the correlation function $R_0$ is invariant under rotation, i.e.\ isotropic, then it as well as the pair correlation function \eqref{eq:gfun} depend only on the spatial distance $\lVert u\rVert$ and the temporal lag $\lvert t\rvert$. The pair correlation function is then given by 
 \begin{align}\label{eq:pcf}
g(u,t)=
 1-\Big\lvert R_0(\|u\|/\alpha_s,\lvert t\rvert/\alpha_t)\Big\rvert^2.
 \end{align}
Further, the space-time $K$-function \citep[see e.g.][]{gabriel:diggle:09,moeller:ghorbani:12} is given by
\begin{align}\label{eq:Kfun}
    K(u,t) &=2\pi\int_{0}^{t}\int_{0}^{u}g(u',t')u'\de u'\de t'\cr
&=\pi u^2 t - \int_{0}^{u}\int_{0}^{t} \lvert R_0(\|u'\|/\alpha_s,\lvert t'\rvert/\alpha_t)\rvert^2 u'\de u'\de t'.
\end{align}

If the intensity $\rho$ is constant, i.e.\ $\rho(u,t)=\rho$, then the process is stationary and the corresponding stationary space-time $g$- and $K$-functions obtain the same formulas \eqref{eq:pcf} and \eqref{eq:Kfun} under the isotropy assumption.
\section{Examples of spatio-temporal determinantal point processes}\label{sec:examples}
Section \ref{sec:cov_and_spectd} recalls the relationship between the covariance function and its spectral density. Then, in Sections~\ref{section:4_1} and~\ref{section:4_2}, the Fourier transform of positive finite measures, i.e. spectral measures, are used to construct stationary spatio-temporal covariance functions, and characteristics of \stDPP{}s with these kernels are derived.

\subsection{The covariance function and its spectral density}\label{sec:cov_and_spectd}

The covariance function of a stationary process can be represented as a Fourier transform of a positive finite measure. 
According to the Wiener-Khintchine theorem \citep[see, e.g.][]{Rasmussen:Williams:06}, if the spectral density function $\varphi(\cdot,\cdot)$ exists, the covariance function $C(\cdot,\cdot)$ and
the spectral density $\varphi(\cdot,\cdot)$ are Fourier duals of each other  given by
\begin{align}\label{eq:daualF}
    C(u,t)&=\int_{\mathbb R^{d+1}}e^{2\pi i (\omega^Tu+ \tau t)} \varphi(\omega,\tau) \de\omega\de\tau,\cr \varphi(\omega,\tau)&=\int_{\mathbb R^{d+1}}e^{-2\pi i(\omega^Tu+\tau t)} C(u,t) \de u\de t,
\end{align}
where $T$ stands for transpose, $\omega$ is the $d$-dimensional spatial component and $\tau$ is  the temporal component. 

A straightforward generalization of  Proposition~3.1 in \cite{Lavancier:etal:15} to the spatio-temporal setting, under continuity and stationarity of $C(u,t)$, when $C(u,t)\in\mathbb{L}^2(\mathbb{R}^d\times \mathbb{R}^+)$,  implies that a \stDPP{} with kernel $C$ exists if the corresponding spectral density satisfies 
\begin{align}\label{eq:prop3.1}
   \varphi(\omega,\tau)\leq 1. 
\end{align}

\subsection{Separable spatio-temporal covariance functions}\label{section:4_1}
A class of separable spatio-temporal covariance functions is usually  given by 
\begin{align}
\label{eq:sepCov}
    C_0(u,t)\propto C_0^s(u)C_0^t(t),
\end{align}
which is valid (i.e.\ a positive definite function) if both the spatial covariance function, $C_0^s(u)$, and the temporal covariance function, $C_0^t(t)$, are valid covariance functions. 
For a separable class of covariance functions, the spectral density $\varphi(\omega,\tau)$ has also a separable form, namely
\begin{align}
\nonumber
    \varphi(\omega,\tau) \propto 
        \Big[\int_{\mathbb{R}^d} e^{-2\pi i\omega^Tu}C_0^s(u) \mathrm{d}u\Big]\times\Big[\int_{\mathbb{R}} e^{-2\pi i\tau t} C_0^t(t) \mathrm{d}t\Big] = \varphi_0^s(\omega)\varphi_0^t(\tau),
\end{align}
where  $\varphi_0^s(\omega)$ and $\varphi_0^t(\tau)$ are the spatial and temporal spectral densities, respectively. According to
 \eqref{eq:prop3.1}, the condition $\varphi_0^s(\omega)\varphi_0^t(\tau) < 1$ must be satisfied for a \stDPP{} with kernel \eqref{eq:sepCov} to exist. 
Further, for this class, the pair correlation function  takes the following simple form 
$$g(u,t)=1-\lvert R_0^s(\|u\|/\alpha_s)\rvert^2\lvert R_0^t(\lvert t\rvert /\alpha_t)\rvert^2,$$
where $R_0^s$ and $R_0^t$ are the correlation functions in space and time  corresponding to  $C_0^s$ and $C_0^t$, respectively.   
Therefore, by \eqref{eq:Kfun} the corresponding $K$-function is given by 
$$
 K(u,t) =\pi u^2 t - \int_{0}^{u} u'\lvert R_0^s(u'/\alpha_s)\rvert ^2 \de u'\int_{0}^{t} \lvert R_0^t(t'/\alpha_t)\rvert ^2\de t'.
$$

There are a large number of classes of valid spatial and valid temporal covariance functions in the literature, for example the Matérn, power exponential and Gaussian classes, to name a few \citep[see, e.g.][]{cressie:wikle:11}. 
As an example, we consider 
the {\em Gaussian covariance function}
$
C_0^s(u)=\sqrt\rho\sigma^2_s\exp(-\|u\|^2/\alpha_s)$, $u\in \mathbb{R}^2$, with spectral density $\varphi_0^s(\omega)=\sqrt\rho\pi\sigma^2_s\alpha_s^2\exp(-\pi^2\alpha_s^2\|\omega\|^2)$,
and the {\em exponential covariance function} 
$
C_0^t(t)=\sqrt\rho\sigma^2_t\exp(-\lvert t\rvert/\alpha_t)
$, 
$
t\in\mathbb{R}^+, 
$
with spectral density $\varphi_0^t(\tau)=(2\sqrt\rho\sigma^2_t \alpha_t)/(1+4\pi^2\alpha_t^2\lvert \tau\rvert^2)$.
Here $\sigma^2_s$ and $\sigma^2_t$ are the variance parameters of the spatial and time components, respectively, and $\alpha_s>0$ and $\alpha_t>0$ are the corresponding range parameters.
A stationary \stDPP{} with intensity $\rho$ and the separable covariance function \eqref{eq:sepCov} with these components, i.e  
\begin{align}
\label{eq:sepCovExam}
    C_0(u,t)=
    \rho\sigma^2_s \sigma^2_t\exp\Big(-\frac{\|u\|^2}{\alpha_s}-\frac{\lvert t\rvert}{\alpha_t}\Big)
\end{align}
will exist if 
\begin{align*}
  \varphi_{sep}(\omega,\tau)=  \frac{2\pi\rho\alpha_s^2\alpha_t\sigma^2_s\sigma^2_t}{(1+4\pi^2\alpha_t^2\tau^2)}\exp{\left(-\pi^2\alpha_s^2\|\omega\|^2\right)}<1.
\end{align*}
Since the maximum of the spectral density occurs at $(0,0)$, so a \stDPP{}($C_0$) exists if $\varphi_{sep}(0,0)<1$, which implies that $\rho<(2\pi\alpha_s^2 \alpha_t\sigma^2_s\sigma^2_t)^{-1}$, and hence the maximal intensity is 
$\rho_{\max}= (2\pi\alpha_s^2 \alpha_t\sigma^2_s\sigma^2_t)^{-1}.$
For this process the pair correlation function is simply given by 
\begin{align}
\label{eq:pcfSep}
 g(u,t)=1-\exp\Big(-\frac{2\|u\|^2}{\alpha_s}-\frac{2\lvert t\rvert}{\alpha_t}\Big).
 \end{align}
The corresponding $K$-function also has a closed-form expression, which is given in \ref{sec:kfun}.
 
\begin{figure}[!ht]
\begin{center}
\includegraphics[width=\textwidth]{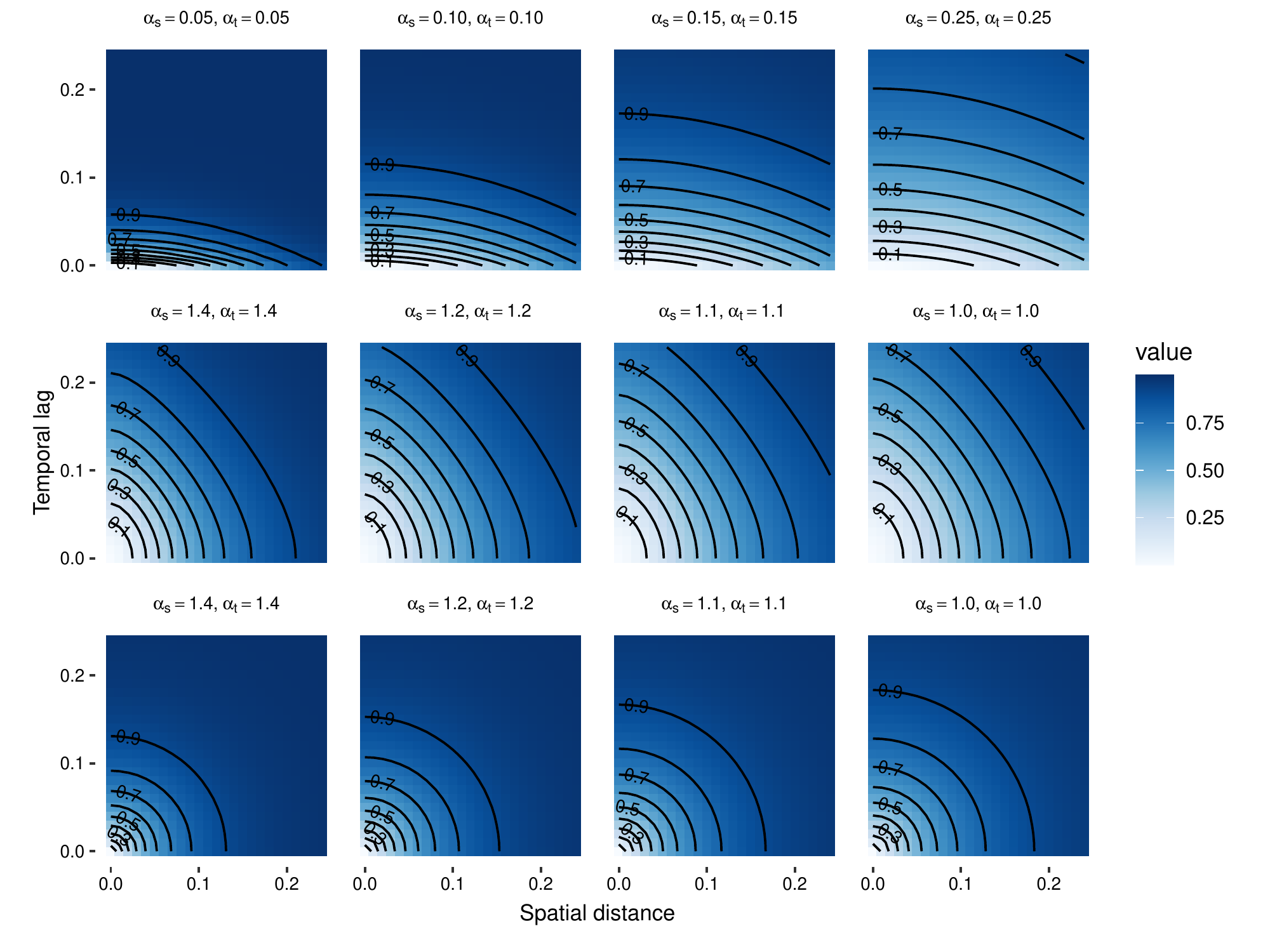}
\caption{Theoretical pair correlation functions  \eqref{eq:pcfSep} (top row), \eqref{eq:pairfunsep} (middle row) and \eqref{eq:pairfun} (bottom row) for different values of the parameters given on top of the plots. }
\label{fig:theopcfs}
\end{center}
\end{figure}

Figure~\ref{fig:theopcfs} (top row) shows that the values of the pair correlation function \eqref{eq:pcfSep} decrease by the increase of the spatial range $\alpha_s$ and the temporal delay $\alpha_t$. Thus, these parameters determine the degree of repulsion for the above separable model.

\subsection{Non-separable spatio-temporal covariance function}\label{section:4_2}
Following \cite{Fuentes:etal:07}, we consider the spatio-temporal spectral density
\begin{align}
    \label{eq:non-sepSpec}
    \varphi_\epsilon(\omega,\tau)=\gamma(\alpha_s^2\alpha_t^2+\alpha_t^2\lvert\omega\rvert^2+\alpha_s^2\tau^2+\epsilon\lvert \omega\rvert^2\tau^2)^{-\nu},
\end{align}
which is an extension  of the commonly used {\em Mat\'{e}rn spectral density} \citep{cressie:wikle:11}.
Here, 
the non-negative parameter $\alpha_s^{-1}$ (spatial range) explains the rate of decay of the spatial correlation, the non-negative parameter $\alpha_t^{-1}$ (temporal delay) 
explains the rate of decay for the temporal correlation. Further, $\gamma>0$ is a scale  parameter.
The parameter $\nu$ measures the degree of smoothness of the process and it should be larger than $(d+1)/2$ to have a well-defined spectral density. The parameter $\epsilon$ controls the interaction between the spatial and temporal components.  For $0\leq\epsilon<1$ the spectral density is non-separable while it is separable when $\epsilon=1$. The maximum of the  spectral density is $\gamma(\alpha_s^2\alpha_t^2)^{-\nu}$ and accordingly  a \stDPP{} with spectral density \eqref{eq:non-sepSpec}  exists if $\gamma <(\alpha_s^2\alpha_t^2)^{\nu}$. 

In the separable case, i.e.\ when $\epsilon=1$, the spectral density is given by
\begin{align}\label{eq:epsilon_one}
   \varphi_{\epsilon=1}(\omega,\tau)=\gamma(\alpha_s^2+\|\omega\|^2)^{-\nu}(\alpha_t^2+\tau^2)^{-\nu} = \varphi_0^s(\omega)\varphi_0^t(\tau).
\end{align}
In this case, $\varphi_0^s(\omega)$ and $\varphi_0^t(\tau)$ are {\em Mat\'{e}rn-type spectral densities} in space and time, respectively. Consequently, the corresponding separable spatio-temporal covariance function, combining \eqref{eq:daualF}, and \eqref{eq:epsilon_one} and using the equations 6.726.4 and 8.432.5 in \cite{Gradshteyn:Ryzhik:07} and setting $d=2$, is given by
\begin{align*}
 C_0(u,t)&=\frac{4\gamma\pi\sqrt\pi}{2^{2\nu-1/2}(\Gamma(\nu))^2}\Big(\frac{2\pi\lvert t\rvert}{\alpha_t}\Big)^{\nu-1/2}\Big(\frac{2\pi\|u\|}{\alpha_s}\Big)^{\nu-1}\cr
 &\times
 \mathcal{K}_{\nu-\frac{1}{2}}\Big(2\pi\alpha_t\lvert t\rvert\Big)\mathcal{K}_{\nu-1}\Big(2\pi\alpha_s\|u\|\Big),
\end{align*}
where $\mathcal{K}_{\nu}(\cdot)$ is the modified Bessel function of the second kind of order $\nu$. $ C_0(u,t)$ is proportional to the product  of  Mat\'{e}rn covariance functions in space and time. 
Using the special cases, $\mathcal{K}_{\frac{1}{2}}(r)=e^{-r}(2r/\pi)^{-1/2}$ and $\mathcal{K}_{\frac{3}{2}}(r)=e^{-r}(1+r^{-1})(2r/\pi)^{-1/2}$ \citep{Abramowitz:Stegun:92}, the above covariance function for $\nu=2$ can be presented as     
\begin{align}\label{kernel}
        C_0(u,t)=\frac{\gamma\pi^2}{4\alpha_s^2\alpha_t^3}(2\pi\alpha_t\lvert t\rvert+1)\exp(-2\pi\alpha_t\lvert t\rvert)(2\pi\alpha_s\|u\|)\mathcal{K}_1(2\pi\alpha_s\|u\|).
\end{align}
For this case, considering the fact that $\lim_{x\rightarrow 0}x\mathcal{K}_1(x)=1$ \citep{Yang:Chu:17}, the intensity of the process is $\rho= C_0(0,0)=(\gamma\pi^2)/(4\alpha_s^2\alpha_t^3)$. Hence, taking into account that $\gamma<\alpha_s^4\alpha_t^4$,  a \stDPP{} with kernel \eqref{kernel} exists if $4\rho\le\pi^2\alpha_s^2\alpha_t $. For this separable case with $\epsilon=1$, considering \eqref{eq:pcf}, the pair correlation function is 
\begin{align}\label{eq:pairfunsep}
    g(u,t)=1- (2\pi\alpha_t\lvert t\rvert+1)^2\exp(-4\pi\alpha_t\lvert t\rvert)\Big[(2\pi\alpha_s\|u\|)\mathcal{K}_1(2\pi\alpha_s\|u\|)\Big]^2.
\end{align} 

For $\epsilon\in(0,1)$ the  spatio-temporal covariance function corresponding to \eqref{eq:non-sepSpec} should be computed numerically as there is no exact closed-form expression. For the case $\epsilon=0$, the stationary non-separable spatial-temporal spectral density is given by
\begin{align}
    \varphi_{\epsilon=0}(\omega,\tau)=\gamma(\alpha_s^2\alpha_t^2+\alpha_t^2\lvert\omega\rvert^2+\alpha_s^2\tau^2)^{-\nu}.
    \label{eq:epsilon_zero}
\end{align} 
Combining \eqref{eq:daualF} and \eqref{eq:epsilon_zero}, and using the equations 6.726.4 and  8.432.5 in \cite{Gradshteyn:Ryzhik:07}, the covariance function when $\epsilon=0$ is 
\begin{align*}
     C_0(u,t)&=\frac{\gamma\pi^{\frac{d+1}{2}}}{2^{\nu-\frac{d+1}{2}}\Gamma(\nu)\alpha_s^{(2\nu-d)}\alpha_t^{(2\nu-1)}}\left\{2\pi\alpha_s\Big((\frac{\alpha_t}{\alpha_s}t)^2+\|u\|^2\Big)^{1/2}\right\}^{\nu-\frac{d+1}{2}}\cr
     &\times \mathcal{K}_{\nu-\frac{d+1}{2}}\left(2\pi\alpha_s\Big((\frac{\alpha_t}{\alpha_s}t)^2+\|u\|^2\Big)^{1/2}\right).
\end{align*}
According to \eqref{eq:prop3.1}, for $d=2$ and $\nu=2$,  there exists a \stDPP{} with kernel
 \begin{align}
 \label{eq:covNonsep}
      C_0(u,t)
 =\frac{\gamma\pi^2}{2\alpha_s^2\alpha_t^3} \exp\left(-2\pi\Big(\alpha_t^2\lvert t\rvert^2+\alpha_s^2\|u\|^2\Big)^{1/2}\right)
 \end{align}
if and only if  
$\gamma<\alpha_s^4\alpha_t^4$. 
Further, it holds that $\rho= C_0(0,0)=(\gamma\pi^2)/(2\alpha_s^2\alpha_t^3)$ for  the covariance functions  \eqref{eq:covNonsep}.
Thus, under the condition $\gamma<\alpha_s^4\alpha_t^4$, it holds that $2\rho\le\pi^2\alpha_s^2\alpha_t$. Therefore, for a \stDPP{} with the above covariance function, the intensity should be at most $\rho_{max}=\pi^2\alpha_s^2\alpha_t/2$.
Further, for this process the pair correlation function is simply given by 
\begin{align}\label{eq:pairfun}
    g(u,t)=1- \exp\left(-4\pi\Big(\alpha_t^2\lvert t\rvert^2+\alpha_s^2\|u\|^2\Big)^{1/2}\right).
\end{align}
The expression for the corresponding $K$-function can be found in \ref{sec:kfun}.

Figure~\ref{fig:theopcfs} (middle and bottom rows) shows that the values of the theoretical pair correlations \eqref{eq:pairfunsep} and \eqref{eq:pairfun} decrease as $\alpha_s^{-1}$ and $\alpha_t^{-1}$ increase. Thus, for these models, the parameters $\alpha_s^{-1}$ and $\alpha_t^{-1}$ play the role of spatial range and time delay that determine the degree of repulsion. 
Moreover, for fixed range parameters, the separable covariance model with \eqref{eq:pairfunsep} leads to smaller values of the pair correlation function and thus more repulsive patterns than the non-separable model with \eqref{eq:pairfun}.  
While the separable covariance function controls to repulsiveness of points in space and time separately, in the non-separable case the points repel each other in the 3D space. This leads to the different small scale interactions.
\section{Discussion and conclusion}\label{sec:discussion}

The different forms of covariance functions presented here allow for \stDPP{}s with different types of repulsion. 
While empirical experiments in the spatio-temporal setting are to be conducted in future work, model fitting for \DPP{}s is available through the maximum likelihood or minimum contrast methods \citep{Lavancier:etal:15} 
based on the summary functions such as the pair correlation function presented here for the given examples, 
and for model assessment, e.g. the global envelope test \citep{myllymaki:etal:17} can be employed.
\section*{Acknowledgement}
The authors are grateful to Frederic Lavancier and Ege Rubak for good discussions. MG was financially supported by the Kempe Foundations (JCSMK22-0134) and 
MM by the Academy of Finland (project numbers 295100 and 327211).

\appendix
\section{$K$-functions of the proposed models}\label{sec:kfun}
Here we give the $K$-functions for two covariance models discussed in Section \ref{sec:examples}. 
The $K$-function of the model with the separable covariance function model \eqref{eq:sepCovExam}, and correspondent to the pair correlation function \eqref{eq:pcfSep}, 
has the following closed-form expression:
\begin{align}\label{eq:Kfun.sep}
    K(u,t) &=\pi u^2 t -\int_{0}^{u}\int_{0}^{t}\exp\Big(\frac{-2u'^2}{\alpha_s}-\frac{2t'}{\alpha_t}\Big)u'\de u'\de t'\cr
    &= \pi u^2 t -\alpha_s\alpha_t/\Big[8\Big(1-e^{(-2u^2/\alpha_s)}\Big)\Big(1-e^{(-2t/\alpha_t)}\Big)\Big].
\end{align}

Employing the general formula of the $K$-function \eqref{eq:Kfun} for the covariance model \eqref{eq:covNonsep}, the space-time $K$-function correspondent to the pair correlation function \eqref{eq:pairfun}  is given by 
\begin{align}\label{eq:Kfun.eps0}
   K(u,t)&= 2\pi\int_{0}^{u}\int_{0}^{t} \Big(1-|R(u'/\alpha_s,t'/\alpha_t)|^2\Big) u'\de u'\de t'\cr
   &= \pi u^2 t -   \frac{2\pi}{\rho^2}\int_{0}^{u}\int_{0}^{t}|C_0(u',t')|^2 u'\de u'\de t'\cr
  &=\pi u^2 t -\frac{\gamma^2\pi^3}{8\rho^2\alpha_s^6\alpha_t^6}\Big[ e^{-4\pi\alpha_t t}\Big(\frac{-2\pi\alpha_t t+e^{4\pi\alpha_t t}-1}{2\pi\alpha_t}\Big)\cr &-\Big(\frac{\alpha_s u}{\alpha_t}+\frac{4\pi\alpha_s^2 u^2}{\alpha_t}\frac{\sinh nt}{2n n!(n+1-(n+1)!)! }\Big)J_1(4\pi,  t)\Big],
\end{align}
where $J_1(4\pi,  t)= \sum_{n=0}^\infty(-2\pi)^n\sum_{m=0}^{n}\frac{\sinh(1+n-2m)t}{m!(n-m!)!(1+n-2m)}$ is an incomplete Bessel function (see more details in \cite{Jones:07}).

\bibliographystyle{dcu} 
\bibliography{paperref.bib}

\end{document}